\documentclass[11pt]{amsart}
\usepackage{amsmath, amsthm, amscd, amssymb}

\newtheorem{thm}{Theorem}[section]
\newtheorem*{thm*}{Theorem}
\newtheorem{cor}[thm]{Corollary}
\newtheorem*{cor*}{Corollary}
\newtheorem{lem}[thm]{Lemma}

\newtheorem*{con*}{Conjecture}
\newtheorem*{prob*}{Problem}

\theoremstyle{definition}

\theoremstyle{remark}
\newtheorem{rem}[thm]{Remark}

\newcommand{\B}{\mathcal{B}}

\begin{document}
\title{On Injective Homomorphisms for Pure Braid Groups,
and Associated Lie Algebras}

\author[F.~R.~Cohen]{F.~R.~Cohen$^{*}$}
\address{Department of Mathematics,
University of Rochester, Rochester, NY 14225}
\email{cohf@math.rochester.edu}
\thanks{$^{*}$Partially supported by the NSF}

\author[Stratos Prassidis]{Stratos Prassidis$^{**}$}
\address{Department of Mathematics
Canisius College, Buffalo, NY 14208, U.S.A.}
\email{prasside@canisius.edu}
\thanks{$^{**}$Partially supported by Canisius College Summer
Grant}
\date{05--11--2004}

\begin{abstract}
The purpose of this article is to record the center of the Lie
algebra obtained from the descending central series of Artin's pure
braid group, a Lie algebra analyzed in work of Kohno
\cite{kohno1,kohno2,kohno3}, and Falk-Randell \cite{falk-rand}. The
structure of this  center gives a Lie algebraic criterion for
testing whether a homomorphism out of the classical pure braid group
is faithful which is analogous to a criterion used to test whether
certain morphisms out of free groups are faithful \cite{cowu1}.
However, it is as unclear whether this criterion for faithfulness
can be applied to any open cases concerning representations of $P_n$
such as the Gassner representation.

\end{abstract}

\maketitle

\section{Introduction}

A classical construction due to Philip Hall dating back to $1933$
gave a Lie algebra associated to any discrete group $\pi$
\cite{hall} which is obtained from filtration quotients of the
descending central series of $\pi$.  That Lie algebra has admitted
applications to the structure of certain discrete groups such as
Burnside groups, as well as applications to problems in topology.
The purpose of this article is to record some additional structure
for this Lie algebra in case $\pi$ is Artin's pure braid group $P_n$
as described below.

That is, define the descending central series of a group $\pi$
inductively by $\{{\Gamma}^k(\pi)\}_{k\ge 1}$ with
\begin{enumerate}
   \item ${\Gamma}^{1}(\pi) = \pi$,
   \item ${\Gamma}^{k}(\pi)$ is the subgroup generated by commutators
   $[ \cdots[\gamma_1,\gamma_2],\gamma_3],\cdots],\gamma_t]$ for $\gamma_i$
   in $\pi$ with $ t \geq k $,
   \item ${\Gamma}^{k+1}(\pi)$ is a normal subgroup of
   ${\Gamma}^{k}(\pi)$,
   \item $E_0^k(\pi) = {\Gamma}^k(\pi)/{\Gamma}^{k+1}(\pi),$ and
   \item $E_0^*(\pi) = \bigoplus_{k\ge 
1}{\Gamma}^k(\pi)/{\Gamma}^{k+1}(\pi).$
\end{enumerate} There is a bilinear homomorphism
$$[-,-]: E_0^p(\pi)\otimes_{\mathbb Z}E_0^q(\pi) \to\ E_0^{p+q}(\pi)$$
induced by the commutator map (not in general a homomorphism) $c:
\pi \times \pi \to\ \pi.$ Natural properties of the map $[-,-]$ due
to P.~Hall, and E.~Witt give $E_0^*(\pi)$ the structure of a Lie
algebra which was developed much further in work of W.~Magnus,
M.~Lazard, A.~I.~Kostrikin, E.~Zelmanov, T.~Kohno \cite{kohno1,
kohno2}, M.~Falk with R.~Randell \cite{falk-rand}, D.~Cohen
\cite{dancohen}, and others.

One standard notation for the Lie algebra attached to the descending
central series is given by $gr_*(\pi)$. The notation $E_0^*(\pi)=
gr_*(\pi)$ used below is adapted from the convention in \cite{
milnor-moore} for the associated graded obtained from a decreasing
filtration.

Let $P_n$ denote the pure braid group on $n$ strands with $B_n$ the
full braid group \cite{magnus-karras-solit,birman}.  A choice of
generators for $P_n$ is $A_{i,j}$, $1 \le i < j \le n$, subject to
relations given in \cite{magnus-karras-solit}. Choices of braids
which represent the $A_{i,j}$ are given by a full twist of strand
$j$ around strand $i$. It is a classical fact using fibrations of
Fadell-Neuwirth (\cite{dane}) that the choice of subgroup generated
by $A_{i,n}$, $1 \le i \le n - 1$, denoted $F_{n-1}$, is free, and
is the kernel of the homomorphism obtained by ``deleting the last
strand".

In the case of the pure braid group, the structure of the Lie
algebra $E_0^*(P_n)$ is given in work of \cite{kohno1},
\cite{kohno2}, and subsequently in \cite{dri}, \cite{falk-rand}:
this Lie algebra is generated by elements $B_{i,j}$ given by the
classes of $A_{i,j}$ in $E_0^1(P_n)$, with $1 \le i < j \le n$.
Since $E_0^1(P_n) = H_1(P_n)$ is an abelian group, the sum of all of
the $B_{i,j}$ given by
$$\Delta(n) = \sum_{1 \le i < j \le n}B_{i,j}$$ is a well-defined element
in $E_0^1(P_n)$. A complete set of relations for $E_0^*(P_n)$, the
``infinitesimal braid relations'', are listed in section $3$ here.

Properties required to state the main result are listed next.
Consider the free group $F[S]$ generated by a set $S$ with $L[S]$
the free Lie algebra generated by the set $S$. A classical fact due
to P.~Hall \cite{hall, serre} is that the morphism of Lie algebras
$e:L[S] \to\ E_0^*(F[S])$ which sends an element $s$ in $S$ to its
equivalence class in $E_0^1(F[S]) = H_1(F[S])$ is an isomorphism of
Lie algebras.

Restrict to the subgroup $F_{n-1}$ the free group generated by
$A_{i,n}$ for $ 1 \leq i < n$. Let $L[V_n]$ denote the free Lie
algebra generated by $B_{i,n}$ with $1 \leq i \le n$. Thus there is
a morphism of Lie algebras
$$\Theta_n: L[V_n] \to\  E_0^*(P_{n})$$ which sends
$B_{i,n}$ to the class of $A_{i,n}$ in $E_0^*(F_{n-1})$. One feature
of $E_0^*(P_{n})$ is that $\Theta_n$ is an isomorphism onto its
image \cite{kohno3,falk-rand}. From now on, $L[V_n]$ is identified
with its image in $E_0^*(P_{n})$.

Let $\mathcal L$ denote a Lie algebra with Lie ideal $\mathcal W$.
The {\it centralizer} of $\mathcal W$  in $\mathcal L$ is defined by
the equation $$C_{\mathcal L}(\mathcal W) = \{x\in \mathcal L|\; [x,
B] = 0, \;\text{for}\; all B \in \mathcal W\}.$$

\begin{thm}\label{thm:centralizers}
If $n > 2$, $$C_{E_0^*(P_n)}(L[V_n]) =
C_{E_0^*(P_n)}(E_0^*(F_{n-1})) = L[{\Delta(n)}].$$
\end{thm}

\begin{rem}
It is quite possible that Theorem \ref{thm:centralizers} appears in
the earlier work concerning the Lie algebra $E_0^*(P_n)$. The
authors are unaware of a reference.
\end{rem}

Several direct corollaries are listed next. Recall the classical
construction of the adjoint representation
$$Ad:L \to Der^{Lie}_*(L)$$ of a graded Lie algebra $L$ for which 
$Der^{Lie}_*(L)$ denotes the
graded Lie algebra of graded derivations of $L$. The map $Ad$ is
defined by the equation $Ad(X)(Y) = [X,Y]$ for $X$, and $Y$ in $L$.
Regard $E_0^*(P_n)$ as a graded Lie algebra by the convention that
$E_0^q(P_n)$ has degree $2q$. Restriction to the Lie ideal $L[V_n]$
gives an induced morphism of Lie algebras $Ad|_{L[V_n]}:E_0^*(P_n)
\to Der^{Lie}_*(L[V_n])$ defined by $Ad|_{L[V_n]}(X)(Y)  = [X,Y].$

\begin{cor}\label{cor:kernel of adjoint representation}
The kernel of the adjoint representation $Ad: E_0^*(P_n)\to
Der_*^{Lie}(E_0^*(P_n))$ as well as the kernel of the restriction of
the adjoint representation $Ad|_{L[V_n]}: E_0^*(P_n)\to
Der_*^{Lie}(L[V_n])$ is given by the cyclic group generated by
$\Delta(n)$ in $E_0^1(P_n)$. Thus there is a short exact sequence of
Lie algebras
\[
\begin{CD}
0 @>{}>> L[\Delta(n)] @>{}>>E_0^*(P_n)@>{Ad|_{L[V_n]}}>>
Image(Ad|_{L[V_n]}) @>{}>> 0.
\end{CD}
\]
\end{cor}

The construction of $E_0^*(\pi)$ is a functor from the category of
discrete groups to the category of Lie algebras over $\mathbb Z$. An
application here, a general method for possibly deciding whether
certain homomorphisms are embeddings, arises from the observation
that if $\pi$ is {\it residually nilpotent}, and a group
homomorphism $f$ out of $\pi$ induces a Lie algebra monomorphism on
the corresponding Lie algebras, then $f$ is a monomorphism
(\cite{cowu1}). This observation is applied to the case of $\pi =
P_n$, the pure braid group on $n$ strands.

\begin{cor}\label{cor:embeddings}
Let $$f: P_n \to\ G$$ be a homomorphism. If the morphisms of Lie
algebras
   $$E_0^*(f)|_{L[V_n]}: L[V_n] \to\  E_0^*(G), \quad
\text{and}\quad E_0^*(f)|_{L[\Delta(n)]}: L[\Delta(n)]\to\
E_0^*(G)$$ are both monomorphisms, then $f$ is a monomorphism. In
addition, the following two statements are equivalent:
\begin{enumerate}
   \item The map $f: P_n \to\ G$ is faithful.
   \item The maps of Lie algebras
   $$E_0^*(f)|_{L[V_n]}: L[V_n] \to\  E_0^*(f(P_n)) \;\;\text{and}\;\;
E_0^*(f)|_{L[\Delta(n)]}: L[\Delta(n)]\to\ E_0^*(f(P_n))$$ are both
monomorphisms where $f(P_n)$ denotes the image of $f$.

\end{enumerate}
\end{cor}

The center $\mathcal Z(n)$ of the braid group $B_n$, isomorphic to
the integers with generator $$(A_{1,2})\cdot (A_{1,3}A_{2,3}) \cdots
(A_{1,n}A_{2,n} \cdots A_{n-1,n})$$
\cite{magnus-karras-solit,birman}, has image in $E_0^*(P_n)$ equal
to $L[\Delta(n)]$. Let
$$i:\mathcal Z(n) \times F_{n-1} \to P_n$$ denote the natural natural 
multiplication map.
The following is a direct consequence of Corollary
\ref{cor:embeddings}.

\begin{cor}\label{cor:centralizer and embeddings}
If the composite
\[
\begin{CD}
  \mathcal Z(n) \times F_{n-1} @>{i}>> P_n @>{f}>> G
\end{CD}
\] induces a monomorphism of Lie algebras

\[
\begin{CD}
L[\Delta(n)] \oplus E_0^*(F_{n-1}) @>{E_0^*(f \circ i)}>> E_0^*(G),
\end{CD}
\] then $f$ is a monomorphism.
\end{cor}

\begin{rem}
Three remarks follow.
\begin{itemize}
     \item A natural question raised by Corollary
\ref{cor:centralizer and embeddings} is as follows. Does the
assumption that $f{\circ}i$ is a monomorphism imply that $f$ is a
monomorphism ? This conclusion does not appear to follow from the
techniques of this paper as the assumption that $f{\circ}i$ is a
monomorphism does not directly imply that $E_0^*(f \circ i)$ is a
monomorphism.
     \item The above Lie algebraic methods for testing whether a 
homomorphism
out of a free group is a monomorphism was used in \cite{cowu1} to
describe certain natural free subgroups of $P_n$ which imply that
the $n$-th homotopy group of the two-sphere is a natural
sub-quotient of $P_n$. \cite{cowu1,bcww}.
     \item It is natural to ask whether the above methods can be
applied to various well-known representations such as the Gassner
representation, the Burau representation for $B_4$ \cite{birman}, or
the Lawrence--Krammer representation (\cite{bigelow},
\cite{krammer}). It is also natural to ask whether there are similar
structure theorems for representations of other related discrete
groups such as those in \cite{dancohen}, or variations which test
whether a representation is both faithful as well as discrete.

The authors have attempted to use the above Lie algebraic methods
above to test whether the classical Burau representation for $B_4$
is faithful. Although there is substantial computer-based evidence
that Corollary \ref{cor:centralizer and embeddings} is satisfied for
the Burau representation of $B_4$, the authors have been unable to
verify this property in general .
\end{itemize}
\end{rem}

  The authors will like to thank Jonathan Pakianathan, Ryan Budney,
and the referee for their suggestions. The referee found a much more
elegant proof for part of the main theorem.

\section{On the Lie algebra for the Pure Braid Group}

The structure of $E_0^*(P_n)$ is given in \cite{kohno1, kohno2,
kohno3}, \cite{falk-rand} and \cite{dri}. Recall that $L[S]$ denotes
the free Lie algebra generated by a set $S$. Then $E_0^*(P_n)$ is
the quotient of the free Lie algebra generated by $B_{i,j}$ for  $1
\leq i < j \leq n$ modulo the {\it infinitesimal braid relations}
(or horizontal 4T relations or Yang-Baxter-Lie relations)
$$E_0^*(P_n) = L[B_{i,j}| \; 1 \le i < j \le n]/I$$
where $I$ denotes the $2$-sided (Lie) ideal generated by the
infinitesimal braid relations as listed next:

\begin{enumerate}
\item $[B_{i,j}, B_{s,t}] = 0$, if $\{i, j\}{\cap}\{s, t\} =
\emptyset$.
\item $[B_{i,j}, B_{i,s} + B_{s,t}] = 0$.
\item $[B_{i,j}, B_{i,t} + B_{j,t}] = 0$.
\item It follows from $2$, and $3$ above that $[B_{j,s}, B_{i,j} +
B_{i,s}] = 0$.
\end{enumerate} In addition, it is convenient to introduce new generators
$B_{j,i}$ for $i<j$ with the convention that
$$B_{j,i} =  B_{i,j}, \;\;\text{for}\;\; i<j.$$

Consider the abelianization homomorphism $$P_n \to P_n/[P_n, P_n] =
E_0^1(P_n) = H_1(P_n)$$ for which the image of $A_{i,j}$ is denoted
$B_{i,j}$. The first homology group $H_1(P_n)$ is isomorphic to
$\oplus_{(n-1)n/2} \mathbb Z$ with basis given by the $B_{i,j}$, for
$1 \leq i < j \leq n$.

Furthermore, there is an induced split short exact sequence of Lie
algebras
$$ 0 \to\ E_0^*(F_{n-1}) \to\ E_0^*(P_n) \to\ E_0^*(P_{n-1}) \to\
0.$$ Thus for each $i>0$, there is a split short exact sequence of
abelian groups $$ 0 \to\ E_0^i(F_{n-1}) \to\ E_0^i(P_n) \to\
E_0^i(P_{n-1}) \to\ 0,$$ and $E_0^i(P_n)$ is isomorphic, as an
abelian group, to $\oplus_{1 \leq j \leq n-1} E_0^i(F_{j})$.

The structure of the Lie algebra $E_0^*(P_n)$ is given in more
detail next via \cite{kohno1,kohno2,kohno3,falk-rand}. Let $L[V_q]$
denote the free Lie algebra ( over $\mathbb Z$ ) generated by the
set $V_q$ with
$$V_q = \{B_{1,q}, B_{2,q}, \cdots, B_{q-1,q}\}, \quad
\text{for}\; 2 \leq q \leq n.$$ Furthermore, there are morphisms of
Lie algebras
$$\Theta_q: L[V_q] \to\ E_0^*(P_n) \quad
\text{given by }\; \Theta_q(B_{j,q}) = B_{j,q}$$ for $1 \leq j < q$
such that the additive extension of the $\Theta_q$ to
$$\Theta: L[V_2] \oplus  L[V_3] \oplus \cdots \oplus
L[V_n] \to\ E_0^*(P_n)$$ is an isomorphism of graded abelian groups.
That is if $a_j(q)$ is an element of $E_0^q(F_{j-1})$ with
$E_0^*(F_{j-1})= L[V_j]$ for $ 2 \leq j \leq n$ and $$x(q) = a_2(q)
+ a_3(q) + \dots + a_{n}(q),$$ then
$$\Theta(x(q)) = \Theta_2(a_2(q)) + \Theta_3(a_3(q)) + \dots +
\Theta_n(a_{n}(q)).$$ The elements $a_j(q)$ will be identified below
with the image $\Theta_j(a_{j}(q))$ unless otherwise noted. The
isomorphism of graded abelian groups $\Theta$ is not an isomorphism
of Lie algebras, but restricts to a morphism of Lie algebras
$\Theta_q: L[V_q] \to\ E_0^*(P_n)$ for each $q \geq 2$. The
infinitesimal braid relations gives the ``twisted'' underlying Lie
algebra structure of $E_0^*(P_n)$.

\begin{lem}\label{lem:ideal}
If $i,j,s <n$, then,
$$[B_{i,j}, B_{s,n}] \in L[B_{1,n}, B_{2,n} \dots B_{n-1,n}] = 
E_0^*(F_{n-1}).$$
Therefore, for each $X \in E_0^*(P_n)$, $[X, B_{s,n}] \in
E_0^*(F_{n-1}),$ and $E_0^*(F_{n-1})$ is a Lie ideal of
$E_0^*(P_n)$.
\end{lem}

\begin{proof}
This follows immediately from the infinitesimal braid relations.
\end{proof}

Centralizers in a free Lie algebra are the subject of the following
exercise from Bourbaki (\cite{bourbaki}, Excercise 3, Chapter II,
section 3 ).

\begin{lem}\label{lem:centralizer}
Let $L[S]$ be the free Lie algebra generated by a set $S$, and let
$a$ be an element of $S$ with $S$ of cardinality at least $2$. Then
the centralizer of $a$ in $L[S]$ is the linear span of $a$.
\end{lem}

\begin{proof}
Let $A_S$ denote the free abelian group generated  by $S$ with $a\in
S$, and $S$ of cardinality at least $2$. The universal enveloping
algebra of $L[S]$ is the tensor algebra $T[A_S]$ while the standard
Lie algebra homomorphism $$j: L[S] \to T[A_S]$$ is injective by the
Poincar\'e-Birkhoff-Witt Theorem (\cite{bourbaki}, and
\cite{jacobson}, p. 168). Identify the elements of $L[S]$ with their
images in $T[A_S]$. Thus if $x\in L[S]$ centralizes $a$, then $a$
commutes with all $x$ in $T[A_S]$.

Consider an element $x$ of $(A_S)^{\otimes n}$ such that  $a \otimes
x = x \otimes a$. Notice that $x = a \otimes x'$ for some element $
x'$ in $(A_S)^{\otimes n - 1}$. Thus by induction on $n$, $x$ is a
scalar multiple of $a^{\otimes n}$, and so $x$ is in the subalgebra
generated by $a$. The intersection of $L[S]$ with the subalgebra
generated by $a$ is precisely the linear span of $a$, thus proving
the lemma.
\end{proof}

The proof of Theorem \ref{thm:centralizers} is given next.

\begin{proof}
There are two parts to this proof. The first part is to show that
the non-zero homogeneous elements of degree $q$ in
$C_{E_0^*(P_n)}(L[V_n])$ are concentrated in degree $q = 1$. The
second part of the proof is to show that the homogeneous elements of
degree $1$ in $C_{E_0^*(P_n)}(L[V_n])$ are precisely scalar
multiples of $\Delta(n) $.

As described above, a restatement of results of Kohno
\cite{kohno1,kohno2}, and Falk-Randell \cite{falk-rand} is that
there is a splitting of $E_0^i(P_n)$ as an abelian group, for each
$i > 0 $:
$$E_0^i(P_n) = E_0^i(L[V_2]) {\oplus}E_0^i(L[V_3]) {\oplus} \dots {\oplus} 
E_0^i(L[V_n])$$
where, for each $1 < m \leq n$, $V_m$ is the linear span of the set
$\{B_{1,m}, B_{2,m}, \dots , B_{m-1,m}\}$.

Let $x(q)$ denote an element in $E_0^q(P_n)$. Thus $x(q)$ is a
linear combination given by $x(q) = a_2(q) + a_3(q) + \dots +
a_n(q)$, $a_j(q) \in E_0^*(L[V_j])$ for which all $a_j(q)$ have the
same degree $q$.

Assume that $x(q)$ is in the centralizer of $L[V_n]$. Thus
$$[x(q), \Gamma ] = 0, \quad \text{for all}\; \Gamma \in L[V_n].$$
It will be shown below by downward induction on $j$ that if $ q >
1$, then $a_j(q) = 0$.

The first case to be checked is that the ``top component'' $a_n(q)$
vanishes for $q > 1$. Assume that $q > 1$. Let $\B(n) = B_{1,n} +
B_{2,n} + \dots + B_{n-1,n}$. The infinitesimal braid relations
$$[B_{i,j}, B_{s,t}] = 0, \quad \text{if}\; \{i, j\}{\cap}\{s, t\} = 
\emptyset$$
and
$$[B_{i,n} + B_{j,n}, B_{i,j}] = 0$$
imply that, for $j < n$, $[a_j(q), \B(n)] = 0$. It follows that
$$[x(q), \B(n)] = [a_n(q), \B(n)] = 0.$$
Thus $a_n(q)$ belongs to the centralizer of the element $\B(n)$ and
both are in $L[V_n]$, which is a free Lie algebra.

By a direct change of basis, there is an equality
$$L[V_n] = L[\B(n), B_{2,n}, \dots , B_{n-1,n}].$$
In addition, notice that Lemma \ref{lem:centralizer} implies that
$a_n(q)$ is a scalar multiple of $\B(n)$ contradicting the
assumption that $q > 1$, and $ n > 2$.

Consider the action of the symmetric group on $n$-letters $\Sigma_n$
on the Lie algebra $E_0^*(P_n)$. This action arises from the
classical action of the symmetric group on $P_n$, and thus induces
automorphisms of the underlying Lie algebra $E_0^*(P_n)$. Note that
this action does not preserve the top free Lie algebra. If $\sigma$
is an element in $\Sigma_n$, then $$ \sigma(B_{i,j}) =
B_{\sigma(i),\sigma(j)} = B_{\sigma(j),\sigma(i)}.$$

By downward induction, assume that $$a_{s+1}(q) = a_{s+2}(q) = \dots
= a_n(q) = 0.$$ Thus $x(q) = a_2(q) + a_3(q) + \dots + a_{s}(q)$ for
$s < n$. Then $$0 = [x(q), B_{s,n}] = [a_{s}(q), B_{s,n}] $$ as
$x(q)$ is assumed to be in the centralizer of $L[V_n]$, and
$[a_i(q),B_{s,n}] = 0$ for $i < s$ by the infinitesimal braid
relations.

Let $\tau_s$ denote the element in $\Sigma_n$ which interchanges
$s$, and $n$ leaving the other points fixed. Regard $\tau_s$ as a
Lie algebra automorphism applied to the previous equation to obtain
$$0 = [\tau_s (a_{s}(q)), \tau_s (B_{s,n})] = [\tau_s (a_{s}(q)),
B_{s,n}].$$ Observe that $\tau_s (a_{s}(q))$ is an element of
$L[V_n]$ as $s < n$, and $\tau_s (a_{s}(q))$ commutes with $B_{s,n}$
with $q > 1$. Hence $\tau_s (a_{s}(q)) = 0$ by Lemma
\ref{lem:centralizer}. Thus, $a_{s}(q) = 0$ as $\tau_s $ is an
automorphism of Lie algebras.

The second part of the proof is an inspection of the homogeneous
elements of degree $1$ in $C_{E_0^*(P_n)}(L[V_n])$, and consists of
showing that these are precisely scalar multiples of $\Delta(n)$ as
is given next. Consider the element $x(1) = a_2(1) + a_3(1) + \cdots
+ a_n(1)$ in $C_{E_0^1(P_n)}(L[V_n])$. Then
$$x(1) = \sum_{1 \leq i < j \leq n} \alpha_{i,j}B_{i,j},$$
where $a_m(1) = \sum_{1 \leq i < m} \alpha_{i,m}B_{i,m}$ for some
choice of integers $\alpha_{i,j}$. Furthermore,
$$[x(1),B_{p,n}]= 0$$ for every $ 1 \leq p < n$ as
$x(1)$ is in $C_{E_0^1(P_n)}(L[V_n])$. It will be checked below that
$\alpha_{i,j} = \alpha_{s,t}$ for all $i<j$, and $s<t$, thus showing
that $x(1)$ is a scalar multiple of $\Delta(n)$.

Notice that $[x(1),B_{p,n}]$ is equal to $$\sum_{i \neq
p,n}\alpha_{i,p}[B_{i,p}, B_{p,n}] + \sum_{i \neq p,n}
\alpha_{i,n}[B_{i,n},B_{p,n}]=  \sum_{i \neq
p,n}(-\alpha_{i,p})[B_{i,n}, B_{p,n}] + \sum_{i \neq p,n}
\alpha_{i,n}[B_{i,n},B_{p,n}]$$ by the infinitesimal braid
relations, and the convention that $B_{i,j} = B_{j,i}$ for $ i < j$.
It follows that $\alpha_{i,n} = \alpha_{i,p}$ as $[B_{i,n},B_{j,n}]$
for $i < j$ form a basis for the homogeneous elements of degree $2$
in $L[V_n]$. A similar computation of $[x(1),B_{p,j}]$ gives
$\alpha_{j,n} = \alpha_{p,j}$ for $p < j <n$.

Thus any element $x(1)$ in the centralizer of $L[V_n]$ is a scalar
multiple of the element $\Delta(n)$. That $\Delta(n)$ centralizes
$E_0^*(P_n)$ follows by inspection. Thus the centralizer of $L[V_n]$
is given by $L[\Delta(n)]$, the free Lie algebra generated by a
single element $\Delta(n)$, a copy of $\mathbb Z$ in degree $1$, and
Theorem \ref{thm:centralizers} follows.
\end{proof}

\section{Proof of Corollaries}

The proof of Corollary \ref{cor:kernel of adjoint representation}
which gives the kernel of the adjoint representation follows next.

\begin{proof}
By definition, the kernel of $Ad|_{L[V_n]}$ is the centralizer of
$L[V_n]$ in $E_0^*(P_n)$, $C_{E_0^*(P_n)}(L[V_n])$. Since
$C_{E_0^*(P_n)}(L[V_n]) = C_{E_0^*(P_n)}(E_0^*(F_{n-1})) =
L[{\Delta(n)}]$ by Theorem \ref{thm:centralizers}, the corollary
follows.
\end{proof}

The next proof is that of \ref{cor:embeddings} which states that if
$f: P_n \to\ G$ is a homomorphism such that the maps of Lie algebras
$$E_0^*(f)|_{L[V_n]}: L[V_n] \to\  E_0^*(G), \;\;\text{and}\;\;
E_0^*(f)|_{L[\Delta(n)]}: L[\Delta(n)]\to\ E_0^*(G)$$ are both
monomorphisms, then $f$ is a monomorphism.

\begin{proof}
Since $P_n$ is residually nilpotent, it suffices to show that
$E_0^*(f)$ is a monomorphism to conclude that $f$ is a monomorphism.

Let $x$ denote an element of least degree $p$ in $E_0^p(P_n)$ which
is in the kernel of $E_0^*(f)$.  Thus $E_0^*(f)([x,B]) = 0$ for
every element $B$ in $E_0^*(P_n)$. Assume that $B$ is in $L[V_n]$.
Then $[x,B]$ is in $L[V_n]$ by Lemma \ref{lem:ideal}. Since
$E_0^*(f)|_{L[V_n]}$ is a monomorphism, $[x,B] = 0$ for every $B$ in
$L[V_n]$. Thus $x$ centralizes $L[V_n]$, and is a multiple of
$\Delta(n)$ by Theorem \ref{thm:centralizers}. Since
$E_0^*(f)|_{L[\Delta(n)]}$ is a monomorphism by hypothesis, $x$ must
be zero.

The final assertion of \ref{cor:embeddings} is the statement that
both $E_0^*(f)|_{L[V_n]}$, and $E_0^*(f)|_{L[\Delta(n)]}$ are
monomorphisms is equivalent to the statement that $f$ is a
monomorphism. This follows directly from the previous step.
\end{proof}

\end{document}